\documentclass[11pt]{article}
\usepackage{amsfonts}
\usepackage{amssymb}
\usepackage{graphicx}
\usepackage{amsmath}
\def\N{{\mathbb N}}

\def\R{{\mathbb R}}

\def\bb{\begin}
\def\bc{\begin{center}}       \def\ec{\end{center}}
\def\be{\begin{equation}}     \def\ee{\end{equation}}
\def\ba{\begin{array}}        \def\ea{\end{array}}
\def\bea{\begin{eqnarray}}    \def\eea{\end{eqnarray}}
\def\beaa{\begin{eqnarray*}}  \def\eeaa{\end{eqnarray*}}
\def\hh{\!\!\!\!}             
\def\EQ{\hh & = & \hh}        
\def\LE{\hh & \le & \hh}

\def\e{\varepsilon}

\def\nn{\nonumber}            \def\ti{\tilde}
\def\oo{\infty}               \def\d{\cdot}
               
\def\lb{\label}               \def\x#1{(\ref{#1})}
\def\q{\quad}                 \def\qq{\qquad}
\def\f{\frac}

\def\B{{\mathcal B}}

\def\M{{\mathcal M}}

\def\Om{\Omega}
\def\OX{\Omega\times X}
\def\ox{\omega, x}
\def\oxk{{\omega_k,x_k}}

\def\dmu{\,{\rm d}\mu}
\def\q{\quad}
\def\qq{\qquad}
\def\qqf{\qquad \forall \ }

\def\BY{\B(Y)}
\def\MX{\M(X)}
\def\MY{\M(Y)}
\def\CX{C(X)}
\def\CY{C(Y)}

\def\Proof{\noindent{\bf Proof} \quad}
\def\qed{\hfill $\Box$ \smallskip}
\def\rd{\,{\rm d}}

\def\bu{$\bullet$\ }

\def\disp{\displaystyle}

\def\lb{\label}
\def\nn{\nonumber}
\def\x#1{(\ref{#1})}

\def\ifl{\iffalse}
\def\Proof{\noindent{\bf Proof} \quad}
\def\qed{\hfill $\Box$ \smallskip}

\def\OX{\Omega\times X}
\def\ox{{\omega,x}}

\begin{document}

\title{\bf On the convergence of the time average for skew-product structure and multiple ergodic system }

\author{\q Xia Pan,$^{1, 2}$  \q Zuohuan Zheng$^{1, 2}$ \q and \q Zhe Zhou$^1$ }

\date{}%

\maketitle

\begin{center}
$^1$ Academy of Mathematics and Systems
Science, Chinese Academy of Sciences, Beijing 100190, China

$^2$ School of Mathematical Sciences, University of Chinese Academy of Sciences, Beijing 100049, China

E-mail: {\tt panxia13@mails.ucas.ac.cn} (Xia Pan) \\
{\tt zhzheng@amt.ac.cn} (Z. Zheng)\\
{\tt zzhou@amss.ac.cn} (Z. Zhou)
\end{center}


\begin{abstract}
In this paper, for a discontinuous skew-product transformation with the integrable observation function, we obtain uniform ergodic theorem and semi-uniform ergodic theorem. The main assumptions are that discontinuity sets of transformation and observation function are neglected in some measure-theoretical sense. The theorems extend the classical results which have been established for continuous dynamical systems or continuous observation functions.  Meanwhile, on the torus $\mathbb{T}^{d}$ with special rotation, we prove the pointwise convergence of multiple ergodic average $\disp \f 1 N \sum_{n=0}^{N-1} f_{1}(R_{\alpha}^{n}x)f_{2}(R_{\alpha}^{2n}x)$  on $\mathbb{T}^{d}$.

\end{abstract}

{\small

{\bf 2010 Mathematics Subject Classification}:
37A20, 
28A35. 

{\bf Key Words and Phrases}: pointwise convergence, uniform ergodicity, semi-uniform ergodicity, multiple ergodic average.
}


\section{Introduction} \lb{intr}
\setcounter{equation}{0}
In 1931, the first major result in ergodic theory was proved by Birkhoff \cite{Bir31} for measure-preserving systems. For topological systems, in 1952, the strongest ergodic theorem, which is called the uniform ergodic theorem, was obtained by Oxtoby \cite{Ox52}. It is a well known result that, for a uniquely ergodic system, time averages of continuous observations converge uniformly. In 1982, a generalization of Oxtoby's ergodic theorem was given by Johnson and Moser \cite{JM82}. By placing severe restrictions on continuous observations instead of the dynamical system, they obtained uniform convergence as well. This was also called the uniform ergodic theorem. As for semi-uniform convergence, in 1997, Stark \cite{S97} gave a semi-uniform version of the sub-additive ergodic theorem for uniquely ergodic systems, though this was not explicitly stated in \cite{S97}. In 2000, Sturman and Stark \cite{SS00} obtained semi-uniform ergodic theorem on a skew-product system. In 2006, Zheng, Xia and Zheng \cite{zx06} gave equivalent conditions for semi-uniform ergodic theorem. Nevertheless, all results mentioned above for topological systems involve continuous dynamical systems with continuous observation functions. In 2011, Zhang and Zhou \cite{ZZ10} extended the uniform ergodic theorem to skew-product quasi-flow, which did not assume the dynamical systems to be continuous in spatial variable. In 2013, Zhang, Zheng and Zhou \cite{ZZ13} established semi-uniform sub-additive ergodic theorem for skew-product transformations, which allows the discontinuity of spatial variable too. Although in \cite{ZZ10,ZZ13} the spatial variable might be discontinuous, observation functions are continuous. In 2011, Dai \cite{Dai11} proved semi-uniform sub-additive ergodic theorem, which allows the discontinuity of observation functions, but spatial variable is continuous.

However, both the spatial variables and observation functions discontinuity, uniform ergodic theorem or semi-uniform ergodic theorem are rare. This paper gives an extension of uniform ergodic theorem and semi-uniform ergodic theorem based on the dynamical system with discontinuous spatial variables, additionally, observation functions are discontinuous.

Next part of this paper concerns pointwise convergence on the torus of the multiple ergodic averages
\be \lb{mul} \f 1 N \sum_{n=0}^{N-1} f_{1}(T_{1}^{n}x)\cdots f_{d}(T_{d}^{n}x).\ee
The convergence of the averages (\ref{mul}) in $L^{2}$ norm was established by Host and Kra \cite{HK05} (see also by Ziegler \cite{Zie07}), where $T_{1}=T,~T_{2}=T^{2},~\cdots,~T_{d}=T^{d}$.  When $T_{1},~T_{2},~\dots,~T_{d}$ are commuting measure-preserving transformations, the convergence of the averages (\ref{mul}) was established by Tao \cite{Tao08}. Soon after, Towsner \cite{Tow09}, Host \cite{BH09} and Austin \cite{Aus10} gave proofs from different viewpoints. Tao's approach was combinatorial and finitary, inspired by the hypergraph regularity and removal lemmas. Towsner used nonstandard analysis, whereas Austin and Host all exploited ergodic methods, building an extension of the original system with good properties. When $T_{1},~T_{2},~\dots,~T_{d}$ belong to nilpotent group, it was proved by Walsh \cite{MW12}. Comparing with $L^{2}$-convergence of the averages (\ref{mul}), the almost surely convergence is relatively few. In 1931, Birkhoff \cite{Bir31} got the almost surely convergence of the averages (\ref{mul}) for $d=1$. The first breakthrough on almost surely convergence of the averages (\ref{mul}) for $d>1$ is due to Bourgain in 1990, who showed in \cite{Bou90} for $d=2$. Subsequently, Huang, Shao and Ye \cite{HSY14} showed the almost surely convergence of the averages (\ref{mul}) for distal systems, Y. Gutman, Huang, Shao and Ye \cite{Ye16} showed the almost surely convergence of the averages (\ref{mul}) pairwise independently determined systems. Recently, E. H. el Abdalaoui \cite{abda17} extended Bourgain's result, he showed that the homogenous ergodic bilinear averages with M\"{o}bius or Liouville weight converged almost surely to zero.

Nevertheless, the results with regards to pointwise convergence of the averages (\ref{mul}) are seldom. Based on Bourgain's theorem in \cite{Bou90}, we give pointwise convergence of the averages (\ref{mul}), this reads as follows:
\be \lb{cir}\f 1 N \sum_{n=0}^{N-1} f_{1}(R_{\alpha}^{n}x)f_{2}(R_{\alpha}^{2n}x)\ee
converges for every point on the torus $\mathbb{T}^{d}$, where $R_{\alpha}=R_{\alpha_{1},\cdots,\alpha_{d}}:\mathbb{T}^{d}\rightarrow \mathbb{T}^{d},$ with $1,\alpha_{1},\cdots,\alpha_{d}$  are rationally independent, $ f_{1}, f_{2}\in C(\mathbb{T}^{d})$.

This paper is organized as follows. We first recall some elements of measure theory and ergodic theory in Section 2. In Section 3, uniform ergodic theorem and semi-uniform ergodic theorem will be given, and detailed proofs will be provided for these theorems. In Section 4, The proof of the pointwise convergence of multiple ergodic average is given.

\section{Preliminary}
\subsection{Existence of invariant measures for SPT}

Let us first recall from \cite{ab17, Wa82} some basic facts on measure theory and ergodic theory. Suppose that $Y$ is a compact metric space. The $\sigma$-algebra of Borel subsets of $Y$ will be denoted by $\B$. All probability measures defined on the measurable space $(Y,\B)$ will be denoted by $\MY$. All real-valued continuous
functions on $Y$ will be denoted by $\CY$. Endowed with the supermum norm $\| \d \|_\oo$,
$(\CY, \| \d \|_\oo)$ is a Banach space. By the Riesz representation theorem \cite[Theorem 6.3]{Wa82},
we know that there exists a bijection between $\MY$ and the set of all normalised positive linear
functionals on $\CY$. Therefore, $\MY$ is identified with a convex subset of the unit ball in $\CY^*$.
Here $\CY^*$ is the dual space of $\CY$. Endowed with the weak$^*$ topology on $\CY^*$, we know that
$(\MY, w^*)$ is convex and sequentially compact.

A measurable map $T: (Y,\B) \to (Y,\B)$ can yield a push $T_*: \MY \to \MY$ by
\be \lb{push}
T_*\mu(B) := \mu(T^{-1}(B))\qqf B\in \B, \ \mu\in \MY. \ee
  Denote the invariant Borel probability measures set by $\M(X,T)=\{\mu\in\MX|T_*\mu=\mu\}$. It is well-known that if $T$ is continuous, then so is the push $T_*$. When $T$ is only measurable, $T_*$ is in general not continuous in weak$^*$ topology. However, when the discontinuous point of $T$ is not too much, we have the following continuity result for $T_*$.

\bb{prop} \lb{pu-con}  {\rm (\cite{Ya04})} Suppose that\/ $T: (Y,\BY) \to (Y,\BY)$ is measurable.
For any sequence\/ $\mu_n \to \mu$ in $(\MY,w^*)$ satisfying \be
\lb{pu-con1} \mu(D_T) =0, \ee one has\/ $T_*\mu_n \to T_* \mu$ in
$(\MY,w^*)$.
\end{prop}

Now, let us introduce a special measurable transformation which admits some invariant measures.

\bb{defn} \lb{SPT} {\rm (\cite{ZZ13})} {\rm Suppose that $\Omega$ and $X$ are compact metric spaces, a transformation\/ $\Phi:\OX \to \OX$
is called a {\it skew-product transformation(SPT)}\/ with a
continuous base\/ $\phi:\Om \to \Om$  if

\bu there holds\/ $\pi \circ \Phi= \phi \circ \pi$ where\/ $\pi: \OX \to \Om$
is the projection and

\bu $\Phi: \OX \to \OX$ is Borel measurable.}
\end{defn}

One of the important features of {\it SPT} is that it allows the discontinuity of $\Phi(\omega, x)$ in the phase space $\OX$. Let
$$D_{\Phi}:=\{ (\ox) \in \OX: \mbox{$\Phi$ is discontinuous at ($\ox$})\}.$$

It is easy to check that $D_{\Phi}$ is a Borel set of $\OX$.
Now combining Proposition \ref{pu-con} with the Schauder-Tychonoff fixed point theorem \cite{DS58},
we obtain the following result which is a complete extension of the classical Bogoliubov-Krylov theorem \cite{N60}.

\bb{prop} \lb{ex-inv}{\rm(\cite{ZZ10})} Let\/ $\Phi$ be an SPT on $\OX$ with a base\/ $\phi$ on\/ $\Om$
fulfilling
$$\nu(\pi(D_{\Phi}))=0 \qqf \nu\in \M(\Om,\phi). \eqno(H)$$
Then there exists at least one invariant Borel probability measure under\/ $\Phi$,
i.e., $\M(\OX,\Phi) \ne \emptyset$. \end{prop}

Some further properties about $\M(\OX,\Phi)$ are listed as follows.

\bb{prop} \lb{fu-inv} {\rm (\cite{ZZ10})} Let\/ $\Phi$ be an SPT on\/ $\OX$ with a base\/ $\phi$ on\/ $\Om$
fulfilling\/ $(H)$. Then

{\rm (i)} $\M(\OX,\Phi)$ is a compact subset of\/ $\M(\OX)$.

{\rm (ii)} $\mu$ is an extreme point of\/ $\M(\OX,\Phi)$ if and only if\/ $\mu$ is ergodic under $\Phi$.
\end{prop}

To prove our main theorems, we need measure$-$theoretic results.

\bb{prop} {\rm(\cite{Bi95})}\lb{con} Let $\mu_{n}$ and $\mu$ be probability measures on $(X,\mathcal{B})$. Then the following two conditions are equivalent:
\begin{enumerate}
  \item $\mu_{n} \rightarrow \mu$ as $n\rightarrow\infty$, in the sense of weak$^*$ topology;
  \item $\lim_{n \to \oo}\int_{X}f(x)\,{\rm d}\mu_{n}(x)=\int_{X}f(x)\,{\rm d}\mu(x)$ for every bounded Borel measurable real function $f(x)$ with $\mu(D_{f})=0$.
\end{enumerate}
\end{prop}

\bb{prop} {\rm(\cite{Wa82})}\lb{invariant} If~ $T:X\rightarrow X$ is continuous and $\mu\in\MX $ then $\mu\in\M(X,T)$ iff $\int f\circ T\,{\rm d}\mu=\int f\,{\rm d}\mu,~~\forall ~f\in \CX$.
\end{prop}

The next results says that each $\mu\in\MX$ is determined by how it integrates bounded and Borel measurable functions. It is a simple extension of classical results for continuous functions, we will just state it here.

\bb{lem} \lb{deng}Let $\mu, \nu$ be two Borel probability measures on the metric space $X$. Then~$\mu=\nu$ iff $\int f\,{\rm d}\mu=\int f\,{\rm d}\nu$, $f$ is bounded and Borel measurable.
\end{lem}

\subsection{Rationally independence rotations of torus}

Consider the unit circle $S^{1}=[0,1]/\sim$, where $\sim$ indicates that $0$ and $1$ are identified. The natural distance on $[0,1]$ induces a distance on $S^{1}$; specifically,
 \[d(x,y)=\min(|x-y|,1-|x-y|).\]
Lebesgue measure on $[0,1]$ gives a natural measure $\mu$ on $S^{1}$, also called Lebesgue measure $\mu$.

  Recall that a measure $\mu$ on the Borel $\sigma$-algebra of a compact topological space $X$ is {\it regular}, if for every $\epsilon>0$ and every $E\in \mathcal{B}(X)$, there is a compact set $M$ and an open set $U$, such that $M\subset E\subset U$ and $\mu(U\backslash M)<\epsilon$. Let $G$ be a compact topological group. There exists a probability measure $\mu$ defined on the Borel $\sigma$-algebra $\mathcal{B}(G)$, such that $\mu(xE)=\mu(E),~~\forall E\in \mathcal{B}(G)$ and $\mu$ is regular. There is only one regular rotation invariant probability measure on $(G,\mathcal{B}(G))$. This unique measure is called {\it Haar measure}. For the circle, the Haar measure is the normalised circular Lebesgue measure. For the torus $\mathbb{T}^{\ell}=\underbrace{S^{1}\times\cdots\times S^{1}}_{\ell~times},~~\ell\geq1$ the Haar measure is the product of the Haar measure on $S^{1}$. 

  A topological dynamical system $f:X\rightarrow X$ is called {\it minimal} if the orbit of every point $x\in X$ is dense in $X$, or, equivalently, if $f$ has no proper closed invariant sets. Let $\ell\geq1$ be an integer and $\alpha=(\alpha_{1},\cdots,\alpha_{\ell})$. The rotation $R_{\alpha}$ has the form
  $$R_{\alpha}(x_{1},\cdots,x_{\ell})=(x_{1}+\alpha_{1},\cdots,x_{\ell}+\alpha_{\ell}) \qq \mbox{(mod 1)}.$$
  For $\ell=1$, and irrational number $\alpha$, the rotation $R_{\alpha}$ is minimal. For $\ell>1$, $\alpha=(\alpha_{1},\cdots,\alpha_{\ell})$, in order to get minimal property of the rotation $R_{\alpha}$, we need auxiliary condition on $\alpha$ as well. 

\bb{defn}{\rm(\cite{Wa82})}
{\rm The real numbers $1,\alpha_{1},\cdots,\alpha_{\ell}$ are {\it rationally independent}, if there is no $k_{0},k_{1},\cdots,k_{\ell}\in Z^{\ell+1}\backslash\{0\}$ such that $k_{0}+k_{1}\alpha_{1}+\cdots k_{\ell}\alpha_{\ell}=0$.}
\end{defn}

\bb{prop} {\rm(\cite{Wa82})}\lb{}  The rotation $R_{\alpha}$ is minimal if and only if the numbers  $\alpha_{1},\cdots,\alpha_{\ell}$ and $1$ are rationally independent.
\end{prop}

From next proposition, the rationally independent rotation $R_{\alpha}=R_{\alpha_{1},\cdots,\alpha_{\ell}}:\mathbb{T}^{\ell}\rightarrow \mathbb{T}^{\ell}$ is uniquely ergodic.
\bb{prop} {\rm(\cite{Wa82})}\lb{unique-ergodicity} Let $\alpha=(\alpha_{1},\cdots,\alpha_{\ell})$ with $1,\alpha_{1},\cdots,\alpha_{\ell}$ rationally independent. The Haar measure is the only probability measure which is invariant by $R_{\alpha}:\mathbb{T}^{\ell}\rightarrow \mathbb{T}^{\ell}.$
\end{prop}

\section{Uniform and semi-uniform convergence for SPT}

In this section $\OX$ will denote a compact metric space. Our aim in this part is to establish some uniform and semi-uniform ergodic theorems for SPT with discontinuous observation functions.

\bb{thm} \lb{equal} Let\/ $\Phi$ be an SPT on\/ $\OX$
with the base\/ $\phi$ on\/ $\Om$ fulfilling\/ $(H)$.
Suppose that\/ $a \in \R$ is a constant, and\/ $f$ is an integrable function with
\be \lb{tj} \mu(D_f)=0 ~~\text{and} ~~  \int_{\OX} f \dmu = a \qqf \mu \in \M(\OX,\Phi).\ee
Then
 \be \lb{} \lim_{n\rightarrow\infty}\disp \f 1 n \sum_{i=0}^{n-1} f ( \Phi^i(\ox))= a \ee
unoformly in $(\ox) \in \OX$.
\end{thm}

 \Proof
  Suppose by contradiction that there is a real number $\epsilon_{0}>0$, a sequence $\{n_{k}\}_{k\geq0}$ of integers tending to $+\infty$ and a sequence $\{(\oxk)\}_{k\geq0}\in \OX$ such that for all $k$,
  \be \lb{non}\left|\f 1 {n_{k}}\sum_{i=0}^{n_{k}-1}f ( \Phi^i(\omega_{k_{j}},x_{k_{j}}))-\int_{\OX} f \dmu\right|\geq\epsilon_{0}.\ee

  Let $\mu_{k}:=\f 1 {n_{k}} \sum_{i=0}^{n_{k}-1}\Phi^{i}_*\delta_{(\oxk)},$ where $\delta_{(\oxk)}$ stands for the Dirac probability measure concentrated at the base point $(\oxk)\in \OX$. By compacity of $\M(\OX)$ in the weak$^*$ topology, one can suppose that the sequence $\{\mu_{k_{j}}\}_{j\geq0}\subseteq\{\mu_{k}\}_{k\geq0}$ converges to a probability measure $\mu'$ which is $\Phi-$invariant.

  We will prove that $\mu'$ is $\Phi-$invariant by the following three steps.

\bu  Step 1: One has
\be\lb{ass1}  \pi_* \mu' \in\M(\Om,\phi).\ee
First of all, we have
\bea
\pi_* \mu'
\EQ \lim_{j\rightarrow\infty}\f 1 {n_{k_{j}}}\sum_{i=0}^{n_{k_{j}}-1}\pi_*(\Phi_*^{i}\delta_{(\omega_{k_{j}},x_{k_{j}})})\nn\\
\EQ \lim_{j\rightarrow\infty}\f 1 {n_{k_{j}}}\sum_{i=0}^{n_{k_{j}}-1}(\pi\circ\Phi^{i})_{*}\delta_{(\omega_{k_{j}},x_{k_{j}})}\nn\\
\EQ \lim_{j\rightarrow\infty}\f 1 {n_{k_{j}}}\sum_{i=0}^{n_{k_{j}}-1}(\phi^{i}\circ\pi)_*\delta_{(\omega_{k_{j}},x_{k_{j}})}.
\nn\eea
In fact, for any $\ti g\in C(\Om)$, define $g=\ti g \circ \pi$. Then we have
\bea \int_{\Om}\ti g\,\rd(\pi_* \mu')
\EQ
\lim_{j\rightarrow\infty}\f 1 {n_{k_{j}}}\sum_{i=0}^{n_{k_{j}}-1}\int_{\Om} \ti g\,{\rm d}(\phi^{i}\circ\pi)_*\delta_{(\omega_{k_{j}},x_{k_{j}})} \nn\\
\EQ
\lim_{j\rightarrow\infty}\f 1 {n_{k_{j}}}\sum_{i=0}^{n_{k_{j}}-1}\int_{\OX} \ti g\circ(\phi^{i}\circ\pi)\,{\rm d}\delta_{(\omega_{k_{j}},x_{k_{j}})} \nn\\
\EQ
\lim_{j\rightarrow\infty}\f 1 {n_{k_{j}}}\sum_{i=0}^{n_{k_{j}}-1}\int_{\OX} \ti g\circ(\pi\circ\Phi^{i})\,{\rm d}\delta_{(\omega_{k_{j}},x_{k_{j}})} \nn\\
\EQ
\lim_{j\rightarrow\infty}\f 1 {n_{k_{j}}}\sum_{i=0}^{n_{k_{j}}-1}\int_{\OX} g\circ\Phi^{i}\,{\rm d}\delta_{(\omega_{k_{j}},x_{k_{j}})}\eea
and
$$\int \ti g\circ \phi\,\rd(\pi_* \mu')=\lim_{j\rightarrow\infty}\f 1 {n_{k_{j}}}\sum_{i=0}^{n_{k_{j}}-1}\int_{\OX} g\circ\Phi^{i+1}\,{\rm d}\delta_{(\omega_{k_{j}},x_{k_{j}})}.$$
Hence

\beaa  \left|\int \ti g\circ \phi\,\rd(\pi_* \mu')-\int \ti g\,\rd(\pi_* \mu')\right|
\EQ \lim_{j\rightarrow\infty}\left|\f 1 {n_{k_{j}}}\int \sum_{i=0}^{n_{k_{j}}-1}( g\circ\Phi^{i+1}-g\circ\Phi^{i})\,{\rm d}\delta_{(\omega_{k_{j}},x_{k_{j}})}\right|\\
\EQ \lim_{j\rightarrow\infty}\left|\f 1 {n_{k_{j}}}\int ( g\circ \Phi^{n_{k_{j}}}-g)\,{\rm d}\delta_{(\omega_{k_{j}},x_{k_{j}})}\right|\\
\LE \lim_{j\rightarrow\infty}\f {2\|g\|_{\infty}} {n_{k_{j}}}\\
\EQ  0.
\eeaa
So $$\int \ti g\circ \phi\,\rd(\pi_* \mu')=\int \ti g\,\rd(\pi_* \mu'),  ~\forall ~\ti g\in C(\Om).$$
According to Proposition \ref{invariant}, we get the desired result \x{ass1}.

\bu Step 2: One has
 \be \lb{ass2} \Phi_* \mu_{k_{j}} \to \Phi_*\mu' \q \mbox{in $(\M(\OX),w^*)$}. \ee
Note that

 \[ 0\le \mu'(D_{\Phi}) \le \mu'(\pi^{-1}(\pi(D_{\Phi})))=
(\pi_*\mu')(\pi(D_{\Phi}))=0, \]
where \x{ass1} and $(H)$ are used. By Propsition \ref{pu-con}, we get the desired result \x{ass2}.
\qed

\bu Step 3: One has
\be \lb{ass3} \Phi_*  \mu_{k_{j}} \to \mu' \q \mbox{in $(\M(\OX), w^*)$}.\ee
In fact, for any $g \in C(\OX)$, we have
\beaa \int_{\OX} g \rd \Phi_* \mu_{k_{j}} \EQ \int_{\OX} g \circ \Phi \rd \mu_{k_{j}}  \\
\EQ \f 1 {n_{k_{j}}}\sum_{i=0}^{n_{k_{j}}-1}\int_{\OX} g\circ\Phi\,{\rm d}\Phi^{i}_*\delta_{(\omega_{k_{j}},x_{k_{j}})}\\
\EQ \f 1 {n_{k_{j}}}\sum_{i=0}^{n_{k_{j}}-1}\int_{\OX} g\circ\Phi\circ\Phi^{i}\,{\rm d}\delta_{(\omega_{k_{j}},x_{k_{j}})}\\
\EQ \disp \f 1 {n_{k_{j}}} \sum_{i=0}^{n_{k_{j}}-1} g \circ \Phi (\Phi^i(\omega_{k_{j}},x_{k_{j}})) \\
\EQ \int_{\OX} g \,{\rm d}\mu_{k_{j}} + \disp \f 1 {n_{k_{j}}} ( g \circ \Phi^{n_{k_{j}}} (\omega_{k_{j}},x_{k_{j}}) - g (\omega_{k_{j}},x_{k_{j}})).
\eeaa
So \[\lim_{j \to \oo}\int_{\OX} g \rd \Phi_* \mu_{k_{j}} = \int_{\OX} g \,{\rm d}\mu'.\]
Since $g \in C(\OX)$ is arbitrary, according to proposition \ref{con}, we get the desired result \x{ass3}.

Combining \x{ass2} with \x{ass3}, we have $\mu' \in \M(\OX, \Phi)$. Finally, it follows from Propsition \ref{con} that \[\lim_{j \to \oo}\int_{\OX} f \rd \mu_{k_{j}}=\int_{\OX}f \,{\rm d}\mu'.\]
Then, we have \[\int_{\OX}f \,{\rm d}\mu'=\lim_{j \to \oo} \f 1 {n_{k_{j}}}\sum_{i=0}^{n_{k_{j}}-1}f ( \Phi^i(\omega_{k_{j}},x_{k_{j}})) \geq a+\e_0,\]
or
\[\int_{\OX}f \,{\rm d}\mu'=\lim_{j \to \oo} \f 1 {n_{k_{j}}}\sum_{i=0}^{n_{k_{j}}-1}f ( \Phi^i(\omega_{k_{j}},x_{k_{j}})) \leq a-\e_0.\]
They are contradiction to condition (\ref{tj}). This ends the proof of Theorem \ref{equal}.
\qed

For a physical problem it is difficult to directly calculate $ \int f \dmu$ for all invariant measures $\mu$. Nevertheless, it is feasible to estimate the range of integration. We will give the semi-uniform ergodic theorem for SPT as follows.

\bb{thm} \lb{inequal1} Let\/ $\Phi$ be an SPT on\/ $\OX$
with the base\/ $\phi$ on\/ $\Om$ fulfilling\/ $(H)$.
Suppose that\/ $a \in \R$ is a constant, and\/ $f$ is an integrable function with
\be \lb{} \mu(D_f)=0 ~~\text{and} ~~ \int_{\OX} f \dmu \leq a \qqf \mu \in \M(\OX,\Phi).\ee
Then, for any given\/ $\e >0$, there exists an $N \in \N$ such that for all\/ $n \geq N$
we have \be \lb{} \disp \f 1 n
\sum_{i=0}^{n-1}  f ( \Phi^i(\ox))\leq a+\e \qqf (\ox) \in \OX. \ee
\end{thm}

Similarly, we have the following result.
\bb{thm} \lb{inequal2} Let\/ $\Phi$ be an SPT on\/ $\OX$
with the base\/ $\phi$ on\/ $\Om$ fulfilling\/ $(H)$.
Suppose that\/ $a \in \R$ is a constant, and\/ $f$ is an integrable function with
\be \lb{} \mu(D_f)=0 ~~\text{and} ~~\int_{\OX} f \dmu \geq a \qqf \mu \in \M(\OX,\Phi).\ee
Then, for any given\/ $\e >0$, there exists an $N \in \N$ such that for all\/ $n \geq N$
we have \be \lb{} \disp \f 1 n
\sum_{i=0}^{n-1}  f ( \Phi^i(\ox))\geq a-\e \qqf (\ox) \in \OX. \ee
\end{thm}

 Repeating the proof of Theorem \ref{equal}, one can prove Theorem \ref{inequal1} and Theorem \ref{inequal2} easily. We omit the details here.

\section{Pointwise convergence on the torus of the multiple ergodic averages} \lb{}
\setcounter{equation}{0}
In this part, we obtain the pointwise convergence with rationally independent rotation on the torus of the multiple ergodic averages.

\bb{thm} \lb{seh} Let $\ell\geq1$ be an integer and $\alpha=(\alpha_{1},\cdots,\alpha_{\ell})$ with $1,\alpha_{1},\cdots,\alpha_{\ell}$ rationally independent, $R_{\alpha}:\mathbb{T}^{\ell}\rightarrow \mathbb{T}^{\ell}$, $f_{1},f_{2}:\mathbb{T}^{\ell}\rightarrow \mathbb{R}$, $f_{1},f_{2}\in C(\mathbb{T}^{\ell})$.
Then $\disp \f 1 N \sum_{n=0}^{N-1} f_{1}(R_{\alpha}^{n}x)f_{2}(R_{\alpha}^{2n}x)$ converges pointwise.
\end{thm}

\Proof
By Bourgain's double recurrence theorem (\cite{Bou90}), for all $f_{1},f_{2}\in C(\mathbb{T}^{\ell})$, the limit $\disp \f 1 N \sum_{n=0}^{N-1} f_{1}(R_{\alpha}^{n}x)f_{2}(R_{\alpha}^{2n}x)$ exists a.e.
 that is to say, there exists
 $$B=\left\{x\in \mathbb{T}^{\ell}\mid \lim_{N \to \oo}\disp \f 1 N \sum_{n=0}^{N-1} f_{1}(R_{\alpha}^{n}x)f_{2}(R_{\alpha}^{2n}x)~~ \mbox{exists}\right\}$$
 with $\mu(B)=1$. Exploiting the method introduced in \cite{Ea17, pan17}, we have $$ \lim_{N \to \oo}\f 1 N \sum_{n=0}^{N-1} f_{1}(R_{\alpha}^{n}x)f_{2}(R_{\alpha}^{2n}x)=\int_{\mathbb{T}^{\ell}} f_{1}\dmu\int_{\mathbb{T}^{\ell}} f_{2}\dmu,~~\forall x\in B,$$
where $\mu$ is the Haar measure of $\mathbb{T}^{\ell}$.

 Let $A=\mathbb{T}^{\ell}\backslash B$. Now, we will assert $A=\emptyset$. Suppose $A\neq\emptyset$.
 Take $y\in A$. Then, there is a real number $\epsilon'>0$, a sequence $(N_{k})_{k\geq0}$ of integers tending to $+\infty$ such that

 \begin{equation}\lb{fan}
  \left|\f 1 {N_{k}}\sum_{n=0}^{N_{k}-1} f_{1}(R_{\alpha}^{n}y)f_{2}(R_{\alpha}^{2n}y)-
 \int_{\mathbb{T}^{\ell}} f_{1}\dmu\int_{\mathbb{T}^{\ell}} f_{2}\dmu\right|\geq\epsilon'.
  \end{equation}
For any $\epsilon>0$, on the one hand, the Haar measure is regular, we can find a $x\in B$, such that $d(x,y)<\epsilon$. On the other hand, the orbit of  any point in $\mathbb{T}^{\ell}$ is dense. There exists an integer $n_{0}>0$, such that $d(R_{\alpha}^{n_{0}}y,x)<\epsilon$. What is more, $R_{\alpha}^{m},~\forall m>0$, is isometric, $d(R_{\alpha}^{n_{0}+m}y,R_{\alpha}^{m}x)=d(R_{\alpha}^{n_{0}}y,x)<\epsilon$. Choose $N_{k}$, such that

 $$\left|\f 1 {N_{k}}\sum_{n=0}^{n_{0}-1} f_{1}(R_{\alpha}^{n}y)f_{2}(R_{\alpha}^{2n}y)\right|+\left|\f 1 {N_{k}}\sum_{n=N_{k}-n_{0}}^{N_{k}-1} f_{1}(R_{\alpha}^{n}x)f_{2}(R_{\alpha}^{2n}x)\right|<\epsilon.$$
 According to the uniformly continuous, there exists a $\delta>0$, if $$d(R_{\alpha}^{n_{0}}y,x)=d(R_{\alpha}^{n_{0}+m}y,R_{\alpha}^{m}x)<\delta,~\forall m>0,$$
 then
  $$d(f_{i}(R_{\alpha}^{n+m}y),f_{i}(R_{\alpha}^{m}x))<\epsilon, ~\text{where}~~ i=1,2.$$
Thus
 \begin{equation}
\begin{aligned}
&\left|\f 1 {N_{k}}\sum_{n=0}^{N_{k}-1} f_{1}(R_{\alpha}^{n}y)f_{2}(R_{\alpha}^{2n}y)-
 \f 1 {N_{k}}\sum_{n=0}^{N_{k}-1} f_{1}(R_{\alpha}^{n}x)f_{2}(R_{\alpha}^{2n}x)\right| \\
&=\left|\f 1 {N_{k}}[\sum_{n=0}^{n_{0}-1} f_{1}(R_{\alpha}^{n}y)f_{2}(R_{\alpha}^{2n}y)+\sum_{n=n_{0}}^{N_{k}-1} f_{1}(R_{\alpha}^{n}y)f_{2}(R_{\alpha}^{2n}y)]\right.\\
&~~~\left.-\f 1 {N_{k}}[\sum_{n=0}^{N_{k}-n_{0}-1} f_{1}(R_{\alpha}^{n}x)f_{2}(R_{\alpha}^{2n}x)+\sum_{n=N_{k}-n_{0}}^{N_{k}-1} f_{1}(R_{\alpha}^{n}x)f_{2}(R_{\alpha}^{2n}x)]\right|\\
&\leq\left|\f 1 {N_{k}}\sum_{n=0}^{n_{0}-1} f_{1}(R_{\alpha}^{n}y)f_{2}(R_{\alpha}^{2n}y)\right|+\left|\f 1 {N_{k}}\sum_{n=N_{k}-n_{0}}^{N_{k}-1} f_{1}(R_{\alpha}^{n}x)f_{2}(R_{\alpha}^{2n}x)\right|\\
&~~~+\f 1 {N_{k}}\left|\sum_{n=n_{0}}^{N_{k}-1} f_{1}(R_{\alpha}^{n}y)f_{2}(R_{\alpha}^{2n}y)-\sum_{n=0}^{N_{k}-n_{0}-1} f_{1}(R_{\alpha}^{n}x)f_{2}(R_{\alpha}^{2n}x)\right|\\
&\leq \epsilon+\f 1 {N_{k}}\left|\sum_{n=0}^{N_{k}-n_{0}-1}[f_{1}(R_{\alpha}^{n+n_{0}}y)f_{2}(R_{\alpha}^{2(n+n_{0})}y)- f_{1}(R_{\alpha}^{n}x)f_{2}(R_{\alpha}^{2n}x)]\right|\\
&=\epsilon+\f 1 {N_{k}}\left|\sum_{n=0}^{N_{k}-n_{0}-1}[(f_{1}(R_{\alpha}^{n+n_{0}}y)-f_{1}(R_{\alpha}^{n}x))f_{2}(R_{\alpha}^{2(n+n_{0})}y)\right.\\
&~~~\left.+ f_{1}(R_{\alpha}^{n}x)(f_{2}(R_{\alpha}^{2(n+n_{0})}y)-f_{2}(R_{\alpha}^{2n}x))]\right|\\
&\leq \epsilon+\epsilon\left|\f 1 {N_{k}}\sum_{n=n_{0}}^{N_{k}-1}f_{2}(R_{\alpha}^{2n}y)\right|+3\epsilon\left|\f 1 {N_{k}}\sum_{n=0}^{N_{k}-n_{0}-1}f_{1}(R_{\alpha}^{n}x)\right|.\\
\end{aligned}
\end{equation}

In fact,
$$\lim_{k \to \oo}\f 1 {N_{k}}\sum_{n=n_{0}}^{N_{k}-1}f_{2}(R_{\alpha}^{2n}y)=\int_{\mathbb{T}^{\ell}} f_{2}\dmu,$$
$$\lim_{k \to \oo}\f 1 {N_{k}}\sum_{n=0}^{N_{k}-n_{0}-1}f_{1}(R_{\alpha}^{n}x)=\int_{\mathbb{T}^{\ell}} f_{1}\dmu.$$

Let
$$C=1+\left|\int_{\mathbb{T}^{\ell}} f_{1}\dmu\right|+3\left|\int_{\mathbb{T}^{\ell}} f_{2}\dmu\right|,$$
$$\epsilon=\f {\epsilon'} C.$$

Then
$$\left|\liminf_{k\rightarrow\infty}\f 1 {N_{k}}\sum_{n=0}^{N_{k}-1} f_{1}(R_{\alpha}^{n}y)f_{2}(R_{\alpha}^{2n}y)-
 \int_{\mathbb{T}^{\ell}} f_{1}\dmu\int_{\mathbb{T}^{\ell}} f_{2}\dmu\right|\leq\epsilon',$$
 and
 $$\left|\limsup_{k\rightarrow\infty}\f 1 {N_{k}}\sum_{n=0}^{N_{k}-1} f_{1}(R_{\alpha}^{n}y)f_{2}(R_{\alpha}^{2n}y)-
 \int_{\mathbb{T}^{\ell}} f_{1}\dmu\int_{\mathbb{T}^{\ell}} f_{2}\dmu\right|\leq\epsilon'.$$
It is contrary to (\ref{fan}). This ends the proof of Theorem \ref{seh}.
\qed

\bb{rem}
If the open problem of almost everywhere convergence for the averages (\ref{mul}) is solved, our methods can be applied to the multiple averages directly, then Theorem {\rm \ref{seh}} holds for all $d\in \N$, i.e. $\disp \f 1 N \sum_{n=0}^{N-1} f_{1}(R_{\alpha}^{n}x)f_{2}(R_{\alpha}^{2n}x)\cdots f_{d}(R_{\alpha}^{dn}x),~~\forall d\in \N$ converges pointwise.
\end{rem}

\noindent{\bf Acknowledgements}: This work is sponsored by 1. National Natural Science Foundation of China(NSFC) under grant 11671382 and 11301512; 2. Key Laboratory of Random Complex Structures and Data Science, Academy of Mathematics and Systems Science, Chinese Academy of Sciences(No.2008DP173182); 3. CAS Key Project of Frontier Sciences(No.QYZDJ-SSW-JSC003); 4. National Center for Mathematics and Interdisciplinary Sciences.


\begin{thebibliography}{222}
\itemsep -2pt

{\small

\bibitem{abda17} E. el Abdalaoui. {\it On the homogeneous ergodic bilinear averages with  M\"{o}bius and Liouville weights}. Preprint, available online at arXiv.org:1706.07280v1 [math.CA].

\bibitem{Ea17} E. el Abdalaoui. {\it On the pointwise convergence of multiple ergodic averages and non-singular dynamical systems}. Preprint, available online at arXiv.org:1406.2608v3[math.DS].

\bibitem{ab17} E. Abdalaoui. {\it The Chowla and the Sarnak conjectures from ergodic theory point of view}. Discrete Contin. Dyn. Syst. 37(2017),2899-2944.

\bibitem{Aus10} T. Austin. {\it On the norm convergence of non-conventional ergodic averages}. Erg. Theory Dynam. Syst., 30(2010),321-338.

\bibitem{Bi95} P. Billingsley. {\it Probability and Measure}.
A Wiley-Interscience Publication, New York, 1995.

\bibitem{Bir31} G. Birkhoff. {\it Proof of the ergodic theorem}. Proc. Natn. Acad. Sci. USA,  17(1931),656-660.

\bibitem{Bou90} J. Bourgain. {\it Double recurrence and almost sure convergence},  J. Reine Angew.Math., 404(1990), 140--161.

\bibitem{Dai11} X. Dai. {\it Optimal state points of the subadditive ergodic theorem}. Nonlinearity, 24(2011),1565--1573.

\bibitem{DS58} N. Dunford and J. Schwartz. {\it Linear Operators, Part I}.
Interscience Publishers, New York, 1958.

\bibitem{Ye16} Y. Gutman, W. Huang, S. Shao and X. Ye. {\it Almost sure convergence of the multiple ergodic averages for certain weakly mixing systems}. Preprint, available online at arXiv.org:1612.02873v1.

\bibitem{HK05} B. Host and B. Kra. {\it Nonconventional ergodic averages and nilmanifolds}. Ann. of Math.,  161(2005),397-488.

\bibitem{BH09} B. Host. {\it Ergodic seminorm for commutting transformations and applications}. Studia Mathematica, 195(2009),31-49.

\bibitem{HSY14} W. Huang, S. Shao and X. Ye. {\it Pointwise convergence of multiple ergodic averages and strictly
ergodic models}. Preprint, available online at arXiv.org:1406.5930.

\bibitem{JM82} R. Johnson and J. Moser. {\it The rotation number
for almost periodic potentials}. Comm. Math. Phys., 84(1982),403-438. Erratum, Comm. Math. Phys., 90(1983),317-318.

\bibitem{Ka95} A. Katok, B. Hasselblatt. {\it Introduction to modern theory of dynamical systems}.
Cambridge University Press, Cambridge/New York, 1995.

\bibitem{N60} V.Nemytskii and V. Stepanov. {\it Qualitative Theory of Differential Equations}.
Princeton:Princeton University Press, 1960.

\bibitem{Ox52} J. C. Oxtoby. {\it Ergodic sets}. Bull. AMS, 58(1952),116-136.

\bibitem{pan17} X. Pan, Z.H. Zheng, Z. Zhou. {\it Ergodic behaviour of nonconventional ergodic averages for commuting transformations}. arXiv:1705.01420 [math.DS].

\bibitem{S97} J. Stark. {\it Invariant graphs for forced systems}. Physica D. 109(1997),163-179.

\bibitem{SS00} R. Sturman and J. Stark. {\it Semi-uniform ergodic theorems
and applications to forced systems}. Nonlinearity, 13(2000),113-143.


\bibitem{Tao08} T. Tao. {\it Norm convergence of multiple ergodic averages for commuting transformations}. Erg. Theory Dynam. Syst., 28(2008),657-688.

\bibitem{Tow09} H. Towsner. {\it Convergence of diagonal ergodic averages}. Erg. Theory Dynam. Syst., 29(2009),1309-1326.

\bibitem{MW12} M. Walsh. {\it Norm convergence of nilpotent ergodic averages}. Ann. of Math. 175(2012),1667-1688.

\bibitem{Wa82} P. Walters. {\it An Introduction to Ergodic Theory},
Springer-Verlag, New York/Berlin, 1982.

\bibitem{Ya04} J. Yan. {\it Lectures on Measure Theory},
Science Press, Beijing, 2004. (in Chinese)

\bibitem{ZZ10} M.R. Zhang and Z. Zhou. {\it Uniform ergodic theorems for discontinuous
skew-product flows and applications to Schr\"{o}dinger equations}. Nonlinearity, 24(2011),1539-1564.

\bibitem{ZZ13} M.R. Zhang, Z.H. Zheng and Z. Zhou. {\it Semi-uniform sub-additive ergodic theorems for discontinuous
skew-product transformations}. Procedings of the American Mathematical Society, 141(2013),3195-3206.

\bibitem{zx06} Z.H. Zheng, J. Xia and Z.M. Zheng. {\it Necessary and suffficent conditions for semi-uniform ergodic theorems and their applications}. Discrete Contin. Dyn. Syst.14(2006),408-417.

\bibitem{Zie07} T. Ziegler. {\it Universal characteristic factors and Furstenberg averages}. J. Amer. Math. Soc., 20(2007),53-97.

}
\end{thebibliography}
\end{document}